\documentclass{article}
\usepackage{amsmath,amssymb}
\usepackage{euscript}
\usepackage{latexsym}
\usepackage[pdftex]{color}
\usepackage[matrix,arrow,curve]{xy}
\usepackage[english]{babel}
\usepackage{wrapfig}
\usepackage{mathrsfs}

\newcommand{\ho}{\underrightarrow{\mathrm{lim}}}

\newtheorem{thm}{\bf Theorem}[section]
\newtheorem{lem}[thm]{\bf Lemma}
\newtheorem{prop}[thm]{\bf Proposition}
\newtheorem{corollary}[thm]{\bf Corollary}
\newtheorem{definition}{\bf Definition}[section]
\newtheorem{example}[definition]{\bf Example}

\def\ge{\geqslant}

\begin{document}

\begin{center}
{\large
	THE HOMOLOGY GROUPS OF RIGHT POINTED SETS OVER A PARTIALLY COMMUTATIVE MONOID
}
\\
\medskip
{
Lopatkin V.
}\\
\end{center}

\begin{abstract}
We study the homology of pointed sets over a partially commutative monoid.
\end{abstract}
2000 Mathematics Subject Classification 18G10, 55U99

\section*{Introduction}
In this paper we consider the homology groups of right pointed sets over a partially commutative monoid $M(E,I)$.\\
\textbf{Definition.} Suppose $E$ is finite set, $I \subseteq E \times E$ is an irreflexive and symmetric relation. A monoid given by a set of generators $E$ and relations $ab = ba$ for all $(a,b) \in I$ is called \emph{free partially commutative} and denoted by $M(E,I)$ \cite{DM}. If $(a,b) \in I$ then members $a,b \in E$ are said to be \emph{commutating generators} \cite{DM}.
\par In this paper we denote by $X^{\bullet}$ a pointed set, in other words $X^{\bullet} = X \cup \{*\}$, here $\{*\}$ is a selected element. The pointed set $X^{\bullet}$ over the partially commutative monoid $M(E,I)$ is also called (see \cite{ol}) $M(E,I)$-set $X^{\bullet}$.

\subsection*{Motivation and bases result}
\par The motivation for the research of homology groups of the $M(E,I)$-sets came from a desire to find topology's invariants for asynchronous transition systems. M. Bednarczyk \cite{B} has introduced \emph{asynchronous transition systems} to the modeling the concurrent processes. In \cite{AT} it was proved that the category of asynchronous transition systems admits an inclusion into the category of pointed sets over free partially commutative monoids. Thus asynchronous transition systems may be considered as $M(E,I)$-sets.
\par It was shown in \cite{HLT} that if the homology of pointed the $M(E,I)$-set are isomorphic to the homology of a point, then the partially commutative monoid $M(E,I)$ consist of one element, and the pointed set consist of one point. Then, we get a following question; suppose that, we have two pointed $M(E,I)$-sets $X^{\bullet}$ and $Y^{\bullet}$, if their homology are isomorphic, does it follow that this $M(E,I)$-sets are isomorphic? In this paper we give the negative answer to this question.

\subsection*{Denoting}
Let $\mathscr{A}$ be a small category. We denote by $\mathrm{Ob} \mathscr{A}$ and $\mathrm{Mor} \mathscr{A}$ the classes of objects and morphisms in $\mathscr{A}$, respectively. By $\mathscr{A}^{op}$ we denote the opposite category.
\par Denote by $\mathrm{Set}_*$ the category of pointed sets: each its objects is a set $X$ with a selected element, written $*$ and called the ''base point``; its morphisms are maps $X \to Y$ which send the base point of $X$ to that of $Y$. Such maps are called \emph{based}. We denote objects of the category $\mathrm{Set}_*$ by $X^{\bullet}$. The set $X^{\bullet} \setminus \{*\}$ denote by $X$.
\par We'll be consider any monoid as the small category with the one object. This exert influence on our terminology. In particular a right $M$-set $X$ will be considered and denoted as a functor $X:M^{op} \to \mathrm{Set}$ (the value of $X$ at the unique object will be denoted by $X(M)$ or shortly $X$.) Morphisms of right $M$-sets are natural transformations. In this paper we consider a case, when pointed set $X^{\bullet}$ is finite. Let us denote by $X_n^{\bullet}$ a finite pointed set that consist of $n+2$ elements $\{x_0,x_1,\ldots,x_n,*\}$. Let us assume that $X_{-1}^{\bullet} = \{*\}$.
\par Suppose that, we have the right pointed $M(E,I)$-set $X_n^{\bullet}$. Let us construct a category $\mathfrak{Cat}_*[M(E,I),X_n^{\bullet}]$, or shortly $\mathfrak{Cat}_*[X_n^{\bullet}]$, which objects are elements of pointed set $X_n^{\bullet}$ and morphisms are triples $(x,\mu,x')$ which $x,x' \in X_n^{\bullet}$ and $\mu \in M(E,I)$ satisfying to $x \cdot \mu = x'$. We assume that the morphism $\left(x,1_{M(E,I)},x  \right)$ is identical morphism of category $\mathfrak{Cat}_*[M(E,I),X_n^{\bullet}]$, here $1_{M(E,I)}$ is identity element of monoid $M(E,I)$. The composition of morphisms $(x,\mu,x')$, $(x',\mu',x'')$ is a morphism of form $(x,\mu \cdot \mu' , x'')$. Sometimes we denote this triples by $x \xrightarrow{\mu} x'$.
\par By $\mathrm{Ab}$ we denote the category of abelian groups and homomorphisms.
\par Let $\mathscr{A}$ be a small category. Denote by $\Delta_{\mathscr{A}}\mathbb{Z}:\mathscr{A} \to \mathrm{Ab}$, or shortly $\Delta \mathbb{Z}$, the functor which has constant values $\Delta \mathbb{Z}(a) = \mathbb{Z}$ at $a \in \mathrm{Ob}\, \mathscr{A}$ and $\Delta \mathbb{Z}(\alpha) = 1_{\mathbb{Z}}$ at $\alpha \in \mathrm{Mor}\,\mathscr{A}$.
\par For any functor $F:\mathscr{A} \to \mathrm{Ab}$, here $\mathscr{A}$ is a small category, denote by $\varinjlim_{k}^{\mathscr{A}} F$ values of the left satellites of the colimit $\varinjlim^{\mathscr{A}}: \mathrm{Ab}^{\mathscr{A}} \to \mathrm{Ab}$. It is well \cite[Prop. 3.3, Aplication 2]{GZ} that there exists an isomorphism of left satellites of the colimit $\varinjlim^{\mathscr{A}}: \mathrm{Ab}^{\mathscr{A}} \to \mathrm{Ab}$ and the functors $H_n(\mathscr{C}_*(\mathscr{A},-)):\mathrm{Ab}^{\mathscr{A}} \to \mathrm{Ab}$. Since the category $\mathrm{Ab}^{\mathscr{A}}$ has enough projectives, these satellites are natural isomorphic to the left derived functor of $\varinjlim^{\mathscr{A}}: \mathrm{Ab}^{\mathscr{A}} \to \mathrm{Ab}$. Denote the values $H_n(\mathscr{C}_*(\mathscr{A},-))$ of satellites at $F \in \mathrm{Ab}^{\mathscr{A}}$ by $\varinjlim_n^{\mathscr{A}}F$.

\section{Basic notions and definitions}\label{bas}
We've assume that, the set of generators of the monoid $M(E,I)$ is finite (see Introduction), hence this monoid $M(E,I)$ satisfying to locally bounded condition \cite{ol}. Then for this monoids there are a following \cite[Theorem 3.1]{ol} 
\begin{thm}\label{general}
Let $M(E,I)$ be a free partially commutative monoid, $\mathcal{X}:M(E,I)^{op} \to \mathrm{Set}$ is a right $M(E,I)$-set, $\mathcal{F}:\mathfrak{Cat}_*[M(E,I),\mathcal{X}] \to \mathrm{Ab}$ is a functor. Then the homology groups $\varinjlim_n^{\mathfrak{Cat}_*[M(E,I),\mathcal{X}]}\mathcal{F}$ are isomorphic to the homology groups of the chain complex
$$0 \xleftarrow{} \bigoplus\limits_{x \in \mathcal{X}}\mathcal{F}(x) \xleftarrow{d_1} \bigoplus\limits_{(x,e_1)}\mathcal{F}(x) \xleftarrow{d_2} \bigoplus_{\substack{(x,e_1e_2),\,e_1 <e_2\\(e_1,e_2) \in I}}\mathcal{F}(x) \xleftarrow{} \ldots$$
$$\ldots \xleftarrow{} \bigoplus_{\substack{(x, e_1,\ldots, e_{n-1}),\,e_1<e_2<\ldots e_{n-1}\\(e_i,e_j)\in I}}\mathcal{F}(x) \xleftarrow{d_n} \bigoplus_{\substack{(x, e_1,\ldots, e_{n}),\,e_1<e_2<\ldots, e_{n}\\(e_i,e_j)\in I}}\mathcal{F}(x) \xleftarrow{} \ldots,$$
with $$d_n(x,e_1,\ldots, e_n,g)=$$
$$=\sum\limits_{s=1}^n(-1)^s\left((x\cdot e_s,e_1,\ldots, \widehat{e_s},\ldots, e_n, \mathcal{F}(x \xrightarrow{e_s}x\cdot e_s)g) - (x,e_1,\ldots, \widehat{e_s},\ldots, e_n, g) \right)$$
\end{thm}
\par From this theorem, in partially, we see that the homology of the $M(E,I)$-set can have a torsion (see \cite[Example 3.2]{ol}). 
\begin{example}\label{EEE}
Let $M(E,I)$ be any free partially commutative monoid with the action over the set $X = \{x_0,*\}$ by the lows
$$x_0 \cdot e = *,\quad * \cdot e = *,$$
for all $e \in E$. Let $\mathbb{Z}[x_0]:\mathfrak{Cat}_*[X_0^{\bullet}] \to \mathrm{Ab}$ be a functor such that $\mathbb{Z}[x_0](x_0) = \mathbb{Z}$ and $\mathbb{Z}[x_0](*) = 0$. By Theorem \ref{general} the groups $\varinjlim_n^{\mathfrak{Cat}_*[X_0^{\bullet}]}\mathbb{Z}[x_0]$ are isomorphic to the homology groups of the chain complex
$$0 \xleftarrow{} \mathbb{Z} \xleftarrow{d_1} \bigoplus_{\substack{e_1}}\mathbb{Z} \xleftarrow{d_2} \bigoplus_{\substack{e_1<e_2\\(e_1,e_2)\in I}}\mathbb{Z} \xleftarrow{d_3} \ldots$$
$$\ldots \xleftarrow{d_{n-1}} \bigoplus_{\substack{e_1<e_2< \ldots < e_{n-1}\\(e_i,e_j)\in I}}\mathbb{Z} \xleftarrow{d_n} \bigoplus_{\substack{e_1<e_2< \ldots < e_{n}\\(e_i,e_j)\in I}}\mathbb{Z} \xleftarrow{d_{n+1}} \ldots,$$
where
$$d_n(e_1,\ldots,e_n) = \sum\limits_{s=1}^n(-1)^{s+1}(e_1,\ldots,\widehat{e_s},\ldots,e_n).$$
Consequently, $\varinjlim_{n}^{\mathfrak{Cat}_*[X_0^{\bullet}]}\mathbb{Z}[x_0] \cong H_{n-1}(E,\mathfrak{M})$ for all $n \ge 2,$ here $H_*(E,\mathfrak{M})$ is homology groups of the simplicial schema $(E,\mathfrak{M})$.
\par In partially, let us assume that, the geometric realization of nonrelation graph of the monoid $M(E,I)$ homeomorphic to the lens space $S^2/\mathbb{Z}_p$, then we get following isomorphism $\varinjlim_{2}^{\mathfrak{Cat}_*[X_0^{\bullet}]}\mathbb{Z}[x_0] \cong \mathbb{Z}_p$.
\end{example}
\par The isomorphism of form
$\ho_n^{\mathfrak{Cat}_*[X_0^{\bullet}]}\mathbb{Z}[x_0] \cong H_{n-1}(E,\mathfrak{M})$, where $n \ge 2$, can be generalized by following
\begin{thm}\label{aug}
Suppose that, we have the pointed set $X_0 = \{x_0,*\}$ over the free partially commutative monoid $M(E,I)$, then there is a following isomorphism for $n \ge 1$
$$\ho_n^{\mathfrak{Cat}_*[X_0^{\bullet}]}\mathbb{Z}[x_0] \cong \widetilde{H}_{n-1}(E,\mathfrak{M}),$$
here $\widetilde{H}_{n-1}(E,\mathfrak{M})$ is reduced homology of the simplicial schema $(E,\mathfrak{M})$.
\end{thm}
\textbf{Proof.} For $n\ge 2$, this isomorphism follows from Example \ref{EEE}. Suppose that, $n \ge 1$. Let us consider a following diagram
$$
  \xymatrix{
\ldots \ar@{->}[r]^(.3){d_3} & \bigoplus\limits_{e_1<e_2,\,(e_1,e_2)\in I}\mathbb{Z} \ar@{->}[r]^(.6){d_2} \ar@{->}[d]_{\cong} & \bigoplus\limits_{e \in E}\mathbb{Z} \ar@{->}[r]^{d_1} \ar@{->}[d]_{\cong} & \mathbb{Z} \ar@{->}[r]^{d_0} \ar@{->}[d]_{\cong} & 0\\
\ldots \ar@{->}[r]^(.3){d_2} & \mathscr{C}_1(E,\mathfrak{M}) \ar@{->}[r]^{d_1} & \mathscr{C}_0(E,\mathfrak{M}) \ar@{->}[r]^(.6){\varepsilon} & \mathbb{Z} \ar@{->}[r] & 0
}
$$
here the top line is the chain complex for the $M(E,I)$-set $X_0^{\bullet}$ with coefficients in the functor $\mathbb{Z}[x_0]$ (see Theorem \ref{general}), and the lower line is augmentation complex for the simplicial schema, and for any vertex $e$ of the $(E,\mathfrak{M})$, $\varepsilon (e) = 1$. Let us remark that, for any $e \in E$, we have $d_1(e) = 1$. Thus, it is not hard to see, that, in dimension $n \ge 0$, the chain complex of the pointed $M(E,I)$-set $X_0$, with coefficients in the functor $\mathbb{Z}[x_0]$, is isomorphic to the augmentation complex $(\mathscr{C}_{n-1}(E,\mathfrak{M}),d_{n-1})$ of the simplicial schema. It means that, there is a following isomorphism, for $n \ge 1$
$$\ho_n^{\mathfrak{Cat}_*[X_0^{\bullet}]}\mathbb{Z}[x_0] \cong \widetilde{H}_{n-1}(E,\mathfrak{M}).$$
\begin{flushright}
Q.E.D.
\end{flushright}
\par From this theorem we get following
\begin{corollary}\label{nottorsion}
The one-dimension homology of the pointed $M(E,I)$-set $X_0^{\bullet}$ with coefficients in the functor $\mathbb{Z}[x_0]$ hasn't the torsion subgroups.
\end{corollary}
\textbf{Proof} is trivial.\\
\par Let us consider the monoid $M(E,I)$ as the small category with the one object. In this case we can to find the homology of the partially commutative monoid $M(E,I)$. Indeed, follow Theorem \ref{general} we get (see \cite[Example 4.1]{ol})
\begin{example}\label{monoid}
Let $M(E,I)$ be a free partially commutative monoid which acts on the one-point set $X_{-1}^{\bullet} = \{X\}$. In this case by Theorem \ref{general} the homology groups $\varinjlim_{s}^{\mathfrak{Cat}_*[X_{-1}^{\bullet}]}\Delta \mathbb{Z}$ are isomorphic to the homology of the chain complex
$$0 \leftarrow \mathbb{Z} \xleftarrow{d_1} \bigoplus_{\substack{e}}\mathbb{Z}  \xleftarrow{d_2} \bigoplus_{\substack{e_1 < e_2\\ (e_1,e_2) \in I}}\mathbb{Z} \leftarrow \ldots$$
$$\ldots \leftarrow \bigoplus_{\substack{e_1< e_2 < \ldots e_{n-1}\\(e_i,e_j) \in I}}\mathbb{Z} \xleftarrow{d_n} \bigoplus_{\substack{e_1< e_2 < \ldots e_{n}\\(e_i,e_j) \in I}}\mathbb{Z} \leftarrow \ldots ,$$
where
$$d_n(e_1,\ldots,e_n) = 0.$$
Consequently, $\varinjlim_{s}^{\mathfrak{Cat}_*[X_{-1}^{\bullet}]}\Delta \mathbb{Z} \cong \mathbb{Z}^{p_s}$, here $p_s$ is the number of subsets $\{e_1,\ldots,e_s\} \subseteq E$ consisting of mutually commutating generators. (The number of empty subsets $p_0 =1$.)
\end{example}

\par Now we'll introduce the following concept. Suppose that, we have the $M(E,I)$-set over the pointed set $X_n^{\bullet}$, we'll say that we have the $M(E,I)$-set over the pointed set $X_n^{\bullet}$ \textbf{with fulling action of the monoid} $M(E,I)$ iff for any $x_p \in X_n^{\bullet}$ and for all $e,e' \in E$, there are following equations $x_p \cdot e = x_p \cdot e'$ (see figure \ref{full}). In addition we'll also say, that the monoid $M(E,I)$ \textbf{acts fully} over the $X_n^{\bullet}$.
\begin{figure}[h!]
$$  
\xymatrix{
&x_0 \ar@/^/@{->}[rd]^(.45){e_1} \ar@/_.8pc/@{->}[rd]_{e_s} \ar@{}[rd]|{{\ldots}} && x_1 \ar@/^/@{->}[ld]^{e_1} \ar@/_.8pc/@{->}[ld]_{e_s} \ar@{}[ld]|{{\ldots}}\\
&& {*} \ar@(r,dr)^(.93){\ldots}^{e_1}  \ar@(l,dl)_{e_s} &
}
\quad
\xymatrix{
&x_0 \ar@/^/@{->}[rd]^(.45){e_2} \ar@/_.8pc/@{->}[rd]_{e_s} \ar@{}[rd]|{{\ldots}} \ar@/^/@{->}[rr]^{e_1} && x_1 \ar@/^/@{->}[ld]^{e_1} \ar@/_.8pc/@{->}[ld]_{e_s} \ar@{}[ld]|{{\ldots}}\\
&& {*} \ar@(r,dr)^(.93){\ldots}^{e_1}  \ar@(l,dl)_{e_s} &
}
$$
\caption{Here is shown actions of the monoid $M(E,I)$ over the pointed set $X_1^{\bullet}= \{x_0,x_1,*\}$; from the left figure we have full action of the monoid, from the right figure the action of the monoid is not fully, because, for example, $x_0 \cdot e_1 \ne x_0 \cdot e_2$.}\label{full1}
\end{figure}
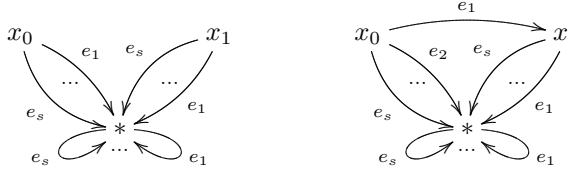

\section{The Homology groups of $M(E,I)$-sets with coefficients in a functor $\mathbb{Z}[x_0,\ldots,x_n]$}\label{full}
For any pointed $M(E,I)$-set $X_n^{\bullet}$, let us introduce a graph $\mathscr{G}[M(E,I),X_n^{\bullet}]$. The set of nodes of a graphs $\mathscr{G}_*[M(E,I),X_n^{\bullet}]$ consists of all elements $\{x_0,\ldots,x_n,*\}$ of the pointed set $X_n^{\bullet}$. For nodes $x_p$, $x_q$ there exists one edge $(x_p,x_q)$ joining $x_p$ to $x_q$, iff for all $e \in E$ there are equations of form $x_p \cdot e = x_q$. A graph $\mathscr{G}[M(E,I),X_n^{\bullet}]$ is a subgraph of the graph $\mathscr{G}_*[M(E,I),X_n^{\bullet}]$ which hasn't the set of edges of form $\xymatrix{{*} \ar@{->}[r]^{(e)} & {*}}$.
\par For any pointed $M(E,I)$-set $X_n^{\bullet}$, let us introduce following conditions \label{con}
\par i) The monoid $M(E,I)$ acts fully over the $X_n^{\bullet}$.
\par ii) The graph $\mathscr{G}[M(E,I),X_n^{\bullet}]$ is a rooting tree with the root in $*$.
\par Suppose we have the partially commutative monoid $M(E,I)$ with full action over the pointed set $X_n^{\bullet} = \{x_0,x_1,\ldots,x_n,*\}$, in other words, we have the $M(E,I)$-set $X_n^{\bullet}$ with full action of the monoid $M(E,I)$. Suppose we know the homology groups of the $M(E,I)$-set $X_0^{\bullet}$ with full action of the monoid. Let us research the homology groups of the $M(E,I)$-set $X_n^{\bullet}$.
\par Let us consider a functor
$$\mathbb{Z}[x_0,\ldots,x_n]: \mathfrak{Cat}_*[\{x_0,\ldots,x_n,*\}] \to \mathrm{Ab}$$
A functor $\mathbb{Z}[x_0,\ldots,x_n]$ is defined on objects by formulas; $\mathbb{Z}[x_0,\ldots,x_n](x_0) = \mathbb{Z}[x_0,\ldots,x_n](x_1) = \ldots = \mathbb{Z}[x_0,\ldots,x_n](x_n) = \mathbb{Z}$ and it is zero in other objects, at morphisms of form $* \xrightarrow{e} *$, this functor is zero, and it is $1_{\mathbb{Z}}$ in other morphisms.
\par Besides, let us introduce a following functor
$$\mathbb{Z}[*]: \mathfrak{Cat}_*[\{x_0,\ldots,x_n,*\}] \to \mathrm{Ab}$$
This functor is $\mathbb{Z}$ on object $*$ and it is zero on other objects, further, on morphisms of form $* \xrightarrow{e} *$, this functor is $1_{\mathbb{Z}}$, and it is zero on other morphisms.  
\par Let us remark that, the category $\mathfrak{Cat}_*[X_{-1}^{\bullet}]$ is full subcategory of the category $\mathfrak{Cat}_*[X_n^{\bullet}]$.
\par Let us prove a following 
\begin{lem}\label{split}
Suppose that, we have the pointed $M(E,I)$-set $X_n^{\bullet}$, then there is a following isomorphism, for $s \ge 1$
$$\ho_s^{\mathfrak{Cat}_*[X_n^{\bullet}]} \Delta_{\mathfrak{Cat}_*[X_n^{\bullet}]} \mathbb{Z} \cong \ho_s^{\mathfrak{Cat}_*[X_n^{\bullet}]} \mathbb{Z}[x_0,\ldots,x_n] \oplus \ho_s^{\mathfrak{Cat}_*[X_{-1}^{\bullet}]} \Delta_{\mathfrak{Cat}_*[X_{-1}^{\bullet}]} \mathbb{Z}.$$
\end{lem}
\textbf{Proof.} From Theorem \ref{general} it isn't hard to see that,
$$\ho_s^{\mathfrak{Cat}_*[X_{-1}^{\bullet}]} \Delta_{\mathfrak{Cat}_*[X_{-1}^{\bullet}]} \mathbb{Z} \cong \ho_s^{\mathfrak{Cat}_*[X_{n}^{\bullet}]} \mathbb{Z}[*] \cong \ho_s^{\mathfrak{Cat}_*[X_{-1}^{\bullet}]}\mathbb{Z}[*].$$
Let us consider a functor $\mathfrak{Ret}:\mathfrak{Cat}_*[X_n^{\bullet}] \to \mathfrak{Cat}_*[X_{-1}^{\bullet}]$. This functor takes all objects of the category $\mathfrak{Cat}_*[X_n^{\bullet}]$ to one object of the category $\mathfrak{Cat}_*[X_{-1}^{\bullet}]$, further, on morphisms of form $(x_p,e,x_q)$, here $x_p \ne *$, this functor is zero, and on other morphisms it is identical. It can easily be checked that the functor $\mathfrak{Ret}$ is left inverse to inclusion functor $\mathfrak{Inj}:\mathfrak{Cat}_*[X_{-1}^{\bullet}] \subseteq \mathfrak{Cat}_*[X_{n}^{\bullet}]$. And since the functor $\mathfrak{Ret}\mathfrak{Inj}: \mathfrak{Cat}_*[X_{-1}^{\bullet}] \to \mathfrak{Cat}_*[X_{-1}^{\bullet}]$ is indicate functor, then the homomorphism $\mathfrak{Ret}_s\mathfrak{Inj}_s$ is indicated automorphism of the group $\ho_s^{\mathfrak{Cat}_*[X_{-1}^{\bullet}]} \mathbb{Z}[*]$, for $s \ge 1$, and moreover there is a following isomorphism
$$\ho_s^{\mathfrak{Cat}_*[X_{n}^{\bullet}]} \Delta \mathbb{Z} \cong \mathrm{Im}\,\mathfrak{Inj}_s \oplus \mathrm{Ker}\,\mathfrak{Ret}_s.$$
This completes the proof.\\
\par Let us show that there is a following
\begin{prop}\label{Z}
Suppose that, we have the pointed $M(E,I)$-set, where the monoid $M(E,I)$ acts fully over the $X_n^{\bullet}$, and assume that the graph $\mathscr{G}[M(E,I),X_n^{\bullet}]$ is the rooting tree with the root in $*$, then there is a following isomorphism, for $s \ge 1$
$$\underrightarrow{\mathrm{lim}}_s^{\mathfrak{Cat}_*[X_n^{\bullet}]}\mathbb{Z}[x_{0},x_{1},\ldots,x_{n}] \cong \left( \underrightarrow{\mathrm{lim}}_s^{\mathfrak{Cat}_*[X_0^{\bullet}]}\mathbb{Z}[x_0] \right)^{n+1}.$$ 
\end{prop}
\textbf{Proof.} Let us consider the category $\mathfrak{Cat}_*[X_n^{\bullet}]$. From conditions i) -- ii) it follows that, there exist $x,\, x' \in \mathrm{Ob}\,\mathfrak{Cat}_*[X_n^{\bullet}]$, and for all $e \in E$, there are a following equations $x \cdot e = x'$, moreover there not exist $x''\in \mathrm{Ob}\,\mathfrak{Cat}_*[X_n^{\bullet}]$, that $x'' \cdot e = x$. Without loss of generality it can be assumed that $x = x_n$, and $x_n \cdot e = x_{n-1}$ for all $e \in E$. Then, it is easily shown that the category $\mathfrak{Cat}_*[\{x_0,\ldots,x_{n-1},*\}] = \mathfrak{Cat}_*[X_{n-1}^{\bullet}]$ is full subcategory of the category $\mathfrak{Cat}_*[X_n^{\bullet}]$. Let us denote by $Q[x_n]$ the inclusion functor $Q[x_n]:\mathfrak{Cat}_*[X_{n-1}^{\bullet}] \subset \mathfrak{Cat}_*[X_{n}^{\bullet}]$. Let us introduce a following functor $P[x_n]:\mathfrak{Cat}_*[X_{n}^{\bullet}] \to \mathfrak{Cat}_*[X_{n-1}^{\bullet}]$. This functor we define in the following way; $P[x_n](x_n) = x_{n-1}$, and on other objects this functor is indicated, further, on morphisms, we assume that $P[x_n](x_n,e,x_{n-1}) = \left(x_{n-1},1_{M(E,I)},x_{n-1} \right)$, on other morphisms this functor is indicated. We see that the functor $P[x_n]$ is left inverse to the inclusion functor $Q[x_n]$. And since the functor $P[x_n]Q[x_n]:\mathfrak{Cat}_*[X_{n-1}^{\bullet}] \to \mathfrak{Cat}_*[X_{n-1}^{\bullet}]$ is indicated functor, thus, for $s \ge 1$, the homomorphism $P_s[x_n]Q_s[x_n]$ is indicated automorphism of the group $\ho_s^{\mathfrak{Cat}_*[X_{n-1}^{\bullet}]}\mathbb{Z}[x_0,\ldots,x_{n}]$, and moreover there is a following isomorphism, for $s \ge 1$
$$\ho_s^{\mathfrak{Cat}_*[X_{n}^{\bullet}]}\mathbb{Z}[x_0,\ldots,x_n] \cong \mathrm{Im}\, Q_s[x_n] \oplus \mathrm{Ker}\,P_s[x_n].$$
It is obvious that $\mathrm{Im}\, Q_s[x_n] \cong \ho_s^{\mathfrak{Cat}_*[X_{n-1}^{\bullet}]}\mathbb{Z}[x_0, \ldots, x_n]$. Also, it is not hard to prove that $\mathrm{Ker}\,P_s[x_n] \cong \ho_s^{\mathfrak{Cat}_*[\{x_n,x_{n-1}\}]}\mathbb{Z}[x_0, \ldots, x_n] \cong \ho_s^{\mathfrak{Cat}_*[X_0^{\bullet}]}\mathbb{Z}[x_0, \ldots, x_n]$. Further, from theorem \ref{general} it follow that $\ho_s^{\mathfrak{Cat}_*[X_{n-1}^{\bullet}]}\mathbb{Z}[x_0, \ldots, x_n] \cong \ho_s^{\mathfrak{Cat}_*[X_{n-1}^{\bullet}]}\mathbb{Z}[x_0, \ldots, x_{n-1}]$ è $\ho_s^{\mathfrak{Cat}_*[X_{0}^{\bullet}]}\mathbb{Z}[x_0, \ldots, x_n] \cong \ho_s^{\mathfrak{Cat}_*[X_{0}^{\bullet}]}\mathbb{Z}[x_0]$. Thus, we have a following isomorphism
$$\ho_s^{\mathfrak{Cat}_*[X_{n}^{\bullet}]}\mathbb{Z}[x_0, \ldots, x_n] \cong \ho_s^{\mathfrak{Cat}_*[X_{n-1}^{\bullet}]}\mathbb{Z}[x_0, \ldots, x_{n-1}] \oplus \ho_s^{\mathfrak{Cat}_*[X_{0}^{\bullet}]}\mathbb{Z}[x_0].$$
In other hand, in the same way, we'll have corresponding isomorphism for categories $\mathfrak{Cat}_*[X_{n-1}^{\bullet}]$, $\mathfrak{Cat}_*[X_{n-2}^{\bullet}]$ e.t.c.  And finally, we'll have isomorphism, for $s \ge 1$
$$\underrightarrow{\mathrm{lim}}_s^{\mathfrak{Cat}_*[X_n^{\bullet}]}\mathbb{Z}[x_{0},x_{1},\ldots,x_{n}] \cong \left( \underrightarrow{\mathrm{lim}}_s^{\mathfrak{Cat}_*[X_0^{\bullet}]}\mathbb{Z}[x_0] \right)^{n+1}$$
\begin{flushright}
Q.E.D.
\end{flushright}
\par The following theorem is our main result of this paper.
\begin{thm}\label{g2}
Suppose, we have the pointed $M(E,I)$-set $X_n^{\bullet} = \{x_0,x_1,\ldots,x_n,*\}$, which satisfaing conditions i) -- ii), then there is a following isomorphism, for  $s \ge 1$
$$H_s(X_n^{\bullet};\mathbb{Z}) \cong 
 \left( \widetilde{H}_{s-1}((E,\mathfrak{M});\mathbb{Z}) \right)^{n+1} \oplus \mathbb{Z}^{p_s},$$
here $p_s$ is the number of subsets $\{e_1,\ldots, e_s\} \subseteq E$ consisting of mutually commutating generators, $\widetilde{H}_{s-1}\left((E,\mathfrak{M}) \right)$ is reduced homology of the simplicial schema $(E,\mathfrak{M})$.
\end{thm}
\textbf{Proof.} Since, conditions i) -- ii) are fulfilled, then from proposition \ref{Z} it follows that, we have a following isomorphism
$$\underrightarrow{\mathrm{lim}}_s^{\mathfrak{Cat}_*[X_n^{\bullet}]} \Delta \mathbb{Z} \cong \left( \underrightarrow{\mathrm{lim}}_s^{\mathfrak{Cat}_*[X_0^{\bullet}]}\mathbb{Z}[x_0] \right)^{n+1} \oplus \underrightarrow{\mathrm{lim}}_s^{\mathfrak{Cat}_*[X_{-1}^{\bullet}]} \Delta \mathbb{Z}.$$
In outher hand, from theorem \ref{aug} it follows that, there is a isomorphism, for $s \ge 1$
$$\underrightarrow{\mathrm{lim}}_s^{\mathfrak{Cat}_*[X_0^{\bullet}]}\mathbb{Z}[x_0] \cong \widetilde{H}_{s-1}(E,\mathfrak{M}).$$ 
But, we already know (see Example \ref{monoid}) that $\ho_s^{\mathfrak{Cat}_*[X_{-1}^{\bullet}]}\Delta \mathbb{Z} \cong \mathbb{Z}^{p_s}$. This completes the proof.\\
\begin{example}
Assume that we have the pointed $M(E,I)$-set $X_3^{\bullet} = \{x_0,x_1,x_2,x_3,*\}$ with fulling action of the monoid $M(E,I)$. This $M(E,I)$-set is shown in figure \ref{fullaction3}.
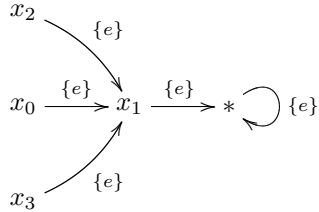
\begin{figure}[h!]
$$
  \xymatrix{
x_2 \ar@/^/@{->}[rd]^{\{e\}} &&&&\\
x_0 \ar@{->}[r]^{\{e\}} & x_1 \ar@{->}[r]^{\{e\}} & {*} \ar@(ur,dr)^{\{e\}}\\
x_3 \ar@/_/@{->}[ur]_{\{e\}} &&&&
}
$$
\caption{Here we sketchy show full action of the monoid $M(E,I)$ over the pointed set $X_3^{\bullet}$. We assume that $e$ runs over all set $E$.}\label{fullaction3}
\end{figure}
\par Let us calculate the homology of this $M(E,I)$-set. Suppose that the geometric realization of the simplicial schema $(E,\mathfrak{M})$ is homeomorphic to a topology space $V$. Thus, from theorem \ref{g2} we get the following isomorphism, for $s \ge 1$
$$H_s(X_{3}^{\bullet}) \cong (H_{s-1}(V))^4 \oplus \mathbb{Z}^{p_s}.$$
\end{example}
\par From Theorem \ref{g2} we get a following
\begin{corollary}
Suppose that, we have the pointed $M(E,I)$-set $X_n^{\bullet}$ which satisfying conditions i) -- ii), then the integer one-dimension homology group of this $M(E,I)$-set hasn't torsion subgroup.
\end{corollary}
\textbf{Proof.} The proof follows from Corollary \ref{nottorsion} and Theorem \ref{g2}.\\
\par From Theorem \ref{g2} follows that, \emph{different} $M(E,I)$-sets can have \emph{isomorphic} homology. Indeed, let us show this by following
\begin{example}\label{differ}
In figure \ref{full1} is shown different $M(E,I)$-sets. We assume that, the set $E$ consist of $s$ elements.
\begin{figure}[h!]
$$
  \xymatrix{
x_0 \ar@/^/@{->}[r]^{e_1} \ar@/_.8pc/@{->}[r]_{e_s} \ar@{}[r]|(.35){\phantom{A^A}^{\vdots}} & x_1 \ar@/^/@{->}[r]^{e_1} \ar@/_.8pc/@{->}[r]_{e_s} \ar@{}[r]|(.35){\phantom{A^A}^{\vdots}} & {*} \ar@(r,ur)^(.65){\vdots}_{e_1}  \ar@(dr,r)_{e_s}
}
\qquad
  \xymatrix{
&x_0 \ar@/^/@{->}[rd]^(.45){e_1} \ar@/_.8pc/@{->}[rd]_{e_s} \ar@{}[rd]|{{\ldots}} && x_1 \ar@/^/@{->}[ld]^{e_1} \ar@/_.8pc/@{->}[ld]_{e_s} \ar@{}[ld]|{{\ldots}}\\
&& {*} \ar@(r,dr)^(.93){\ldots}^{e_1}  \ar@(l,dl)_{e_s} &
}
$$
\caption{Here is shown different actions of the monoid $M(E,I)$ over the pointed set $X_1^{\bullet} = \{x_0,x_1,*\}$.}\label{full1}
\end{figure}
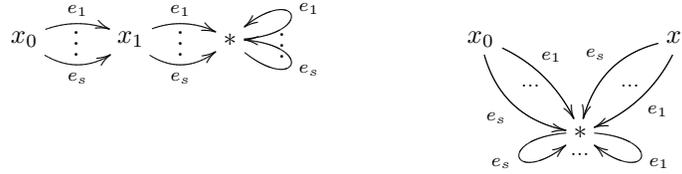
But this different $M(E,I)$-sets have isomorphic homology. Indeed, let us denote the $M(E,I)$-set from left figure by $\mathfrak{L}$, and the $M(E,I)$-set from right figure, we denote by $\mathfrak{R}$. Then, from Theorem \ref{g2}, we get the following isomorphisms, for $n \ge 1$
$$H_n(\mathfrak{L};\mathbb{Z}) \cong \left( \widetilde{H}_{n-1}(E,\mathfrak{M}) \right)^2 \oplus \mathbb{Z}^{p_n};$$
$$H_n(\mathfrak{R};\mathbb{Z}) \cong \left( \widetilde{H}_{n-1}(E,\mathfrak{M}) \right)^2 \oplus \mathbb{Z}^{p_n};$$
The null-dimension homology also isomorphic as $M(E,I)$-sets have common connected components.
\end{example}

\section*{Concluding remarks}
We've see that, different $M(E,I)$-sets can have isomorphic homology groups. This means that, homology groups is not enough for research of topology's invariants of $M(E,I)$-sets (asynchronous transition systems). Thus there arise a motivation for research of others objects of algebraic topology for pointed $M(E,I)$-sets.

\begin{flushleft}
Lopatkin Viktor \quad wickktor@gmail.com\\
Komsomolsk on Amure State Technical University\\
Russia
\end{flushleft}

\end{document}